\documentclass[12pt]{amsart}
\usepackage{amssymb,amsmath,eufrak,amscd}
\usepackage[mathscr]{eucal}

\theoremstyle{plain}
\newtheorem{theorem}{Theorem}[section]
\newtheorem{proposition}[theorem]{Proposition}
\newtheorem{corollary}[theorem]{Corollary}
\newtheorem{lemma}[theorem]{Lemma}

\newtheorem{proposition.a}[theoremaa]{Proposition}
\newtheorem{corollary.a}[theoremaa]{Theorem}
\newtheorem{ccorollary.a}[theoremaa]{Corollary}

\theoremstyle{definition}
\newtheorem{definition}[theorem]{Definition}
\newtheorem{remark}[theorem]{Remark}
\newtheorem{example}[theorem]{Example}

\newtheorem{discussion}[theorem]{Discussion}

\newtheorem{convention}[theorem]{Convention}

\newcommand{\ra}{\rightarrow}
\newcommand{\lra}{\longrightarrow}
\newcommand{\noi}{\noindent}
\newcommand{\PP}{\mathbf{P}}

\newcommand{\RR}{\mathbf{R}}
\newcommand{\NN}{\mathbf{N}}
\newcommand{\A}{\mathbf{A}}
\newcommand{\ZZ}{\mathbf{Z}}
\newcommand{\CC}{\mathbf{C}}
\newcommand{\QQ}{\mathbf{Q}}

\newcommand{\OO}{\mathcal  {O}}

\newcommand{\II}{\mathcal  {I}}

\newcommand{\Spec}{\textnormal{Spec}}
\newcommand{\fra}{\frak{a}}
\newcommand{\frb}{\frak{b}}

\newcommand{\frm}{\frak{m}}
\newcommand{\frmm}{\frak{m}}
\newcommand{\frakm}{\frmm}
\newcommand{\frc}{\frak{c}}

\newcommand{\frj}{\frak{j}}

\newcommand{\frd}{\frak{d}}

\newcommand{\hh}[3]{h^{{#1}} \big( {#2} , {#3} \big) }

\newcommand{\for}{ \ \ \text{ for } \ }
\newcommand{\fall}{ \ \ \text{ for all } \ }

\newcommand{\length}{\textnormal{length}}
\newcommand{\vol}{\textnormal{vol}}

\newcommand{\MI}[1]{\mathcal  {J} ( {#1} ) }

\newcommand{\pr}{\prime}

\newcommand{\spec}{\text{Spec}}

\title[Approximation of Valuation Ideals]{Uniform
Approximation of Abhyankar Valuation Ideals in Smooth
Function Fields} 

\author{Lawrence Ein}
\address{Department of Mathematics \\University of
Illinois at Chicago \hfil\break\indent  851 South Morgan
Street (M/C 249)\\ Chicago, IL 60607-7045, USA}
\email{ein@math.uic.edu}
\thanks{Research of Ein partially supported by NSF Grant
DMS  99-70295}

\author{Robert Lazarsfeld} 
\address{Department of Mathematics
\\ University of Michigan \\ Ann Arbor, MI 48109, USA}
\email{rlaz@math.lsa.umich.edu}
\thanks{Research of  Lazarsfeld  partially supported
by   NSF Grant DMS 97-13149}

\author{Karen E. Smith}
\address{Department of Mathematics
\\ University of Michigan \\ Ann Arbor, MI 48109, USA}
\email{kesmith@math.lsa.umich.edu}
\thanks{Research of Smith partially supported by  NSF
Grant DMS 00-70722, and by a US Fulbright Fellowship}

\begin{document}

\maketitle


\section*{Introduction}

In this paper we use the theory of multiplier ideals
 to  show that the valuation ideals of a rank one
Abhyankar  valuation centered at a smooth point of a
complex algebraic variety 
 are approximated, in a quite strong sense,  by 
 sequences of powers of fixed ideals.  

Let  $R$ be an
$n$-dimensional  regular local domain essentially of
finite type over a ground field
$k$ of characteristic zero, and let
$\nu$ be a rank one  valuation centered on
$R$. Recall that this is equivalent to asking that $\nu$
be an
$\RR$-valued valuation on the fraction field $K$ of $R$,
taking non-negative values on $R$ and positive values on
the maximal ideal
$\frakm
\subseteq R$. A theorem of Zariski and Abhyankar states
that  
\begin{equation} \label{Zariski-Abhyank.Ineq}
{\text{trans.deg }} \nu \, + \,  {\text{rat.rank
}} \nu \ \le \ \dim K/k,
\end{equation}
where the \textit{rational rank} of $\nu$ is the rank of
its  value group, while its \textit{transcendence degree}
 is the maximal dimension of the center
$\nu$  on some model of $K / k$. One says that $\nu$ is
an \textit{Abhyankar valuation} if equality holds in
(\ref{Zariski-Abhyank.Ineq}), i.e. if
\[
{\text{trans.deg }} \nu \, + \,  {\text{rat.rank
}} \nu \ = \ \dim K/k.
\]
Among all the valuations centered on $R$, these are the
ones from which one expects the best behavior. 
 For
example, a  divisorial valuation --- corresponding to
the order of vanishing along a prime divisor
$E \subseteq Y$ contracting to the closed point of
$X = \Spec(R)$ under a proper birational map $Y \lra X$
--- is an Abhyankar valuation centered on $R$ having transcendence
degree $n -1$ and rational rank $1$. At the other
extreme, if $\alpha_1, \ldots , \alpha_n \in \RR$ are
any $\QQ$-linearly independent positive real numbers,
then there exists a unique valuation $\nu$ centered on
$k[x_1, \ldots, x_n]_{(x_1,\ldots, x_n)}$ with $\nu(x_i)
= \alpha_i$: here  $\text{trans.deg } \nu = 0$ and 
$\text{rat.rank }  \nu = n$. 
Abhyankar valuations have been the focus of considerable
attention. For example, they are known to 
admit local uniformization in any characteristic \cite{KK}, and 
 already when
$\dim R = 2$ they involve a great deal of beautiful and
intricate geometry \cite{Spiv}.

Given a valuation  $\nu$ as above, denote by
$\Phi =
\nu(R)
\subseteq \RR$ the value semigroup of $\nu$ on $R$. For
each real number
$m \in
\Phi$, 
 let 
\[ \fra_m  \ = \ \big \{ f \in R \mid \nu(f) \ge m
\big \}
\]
 denote the ideal of $R$ consisting of all elements of
$R$ whose values  are at least $m$. Clearly $\fra_m^\ell
\subseteq \fra_{\ell m}$ for every natural number $\ell
\in \NN$, but typically the inclusion is strict. However
our main theorem shows that for Abhyankar valuations
these two ideals lie surprisingly close to each other: 
\begin{corollary.a} \label{New.Main.Thm} Let $\nu$ be an Abhyankar
valuation centered on
$R$. Then  there exists a  fixed value $e \in \Phi$ such that 
$$
\fra_m^{\ell} \ \subseteq \ \fra_{m\ell}\  \subseteq \  \fra_{m-e}^{\ell}
$$
for all $m \in \Phi$ and all $\ell \in \NN$.\footnote{We adopt the
convention that $\fra_{m-e} = R$ when $m < e$.} 
\end{corollary.a}
\noi Roughly speaking, the Theorem asserts that  the valuation ideals
$\fra_{m
\ell}$ are closely and uniformly approximated by powers of $\fra_m$.
It follows from the Theorem that there is a fixed non-zero element $\delta
\in R$ such that
\begin{equation}\label{cor}
\delta^{\ell} \fra_{m\ell} \ \subseteq \  \fra_m^{\ell} 
\end{equation}  for all $m \in \Phi$ and  
$\ell \in \NN$. This points to the heuristic idea that the associated
``Rees ring" of the Abhyankar valuation
$\nu $,  while   usually not finitely generated
when $n = \dim R > 1$, is 
``almost" finitely generated. It would be interesting to
know if one can make this precise. Both the Theorem and 
(\ref{cor}) can
fail  for non-Abhyankar valuations:  see Remark \ref{fail}.

Theorem A can be interpreted as a  strengthened
form of a celebrated theorem of Izumi.  In the
setting of Theorem  A,
 Izumi's result (in a form due to H\"ubl and
Swanson)
 asserts that there
exists an index
$p$ such that
$\fra_{p\ell}
 \subseteq  
\frmm^{\ell}$ for all natural numbers $\ell$,
  where as above $\frmm$ is the maximal ideal of
$R$.\footnote{Izumi's statement actually  deals with less
general valuations  but  more general rings. It is
proved in 
\cite{Iz},  \cite{Rees}, and
\cite{Hub.Swan}.} 
Theorem \ref{New.Main.Thm} says
much more:
 not only is $\fra_{p\ell} \subseteq \frmm^{\ell}$ for
all
$\ell$, but  $\fra_{p\ell}$ is contained in the
$\ell$-th power of an ideal $\fra_{p-e}$ very close to 
$\fra_p$. In other words
$\fra_{p\ell} \subseteq \fra_{p-e}^{\ell}$,  where not only
are the 
 $\fra_{p-e}$'s getting deeper in $R$  as a function of $p$ (whereas
$\frmm$ stays fixed),   but  they are getting deeper at the same rate
as the $\fra_p$ themselves. 
The theorem also implies---at least for regular local rings
essentially of finite type over a field--- the well-known statement of 
Izumi's Theorem comparing values of two valuations centered at $\frmm$;
see Corollary \ref{izumi}.

In an algebro-geometric context, Theorem A is most interesting 
for divisorial valuations. More generally, Theorem A can be applied to
a composite of valuations to yield the following:
\begin{ccorollary.a}
Let $D$ be an effective divisor on a normal variety $X$ and suppose
that $X \overset{\pi}\to Y$ is a proper birational map contracting $D$ to a
smooth (but not necessarily closed) point of $Y$.
Then there exists a natural number $e$ such that 
$$
\pi_*\OO_X(-m\ell D) \subset [\pi_*\OO_X(-(m-e)D)]^{\ell}
$$
for all natural numbers $\ell$ and all $m \geq e$.
\end{ccorollary.a}
\noi
Thus even though the algebra 
\begin{equation}\label{alg}
 \bigoplus_{\ell \in \NN} \pi_*\OO_X(-m\ell D)
\end{equation}
is not finitely generated in general, Corollary B gives some measure
of control over it. 
 Of course, it is a central problem of birational
geometry to understand when such algebras are finitely generated:
when (\ref{alg}) is  finitely generated,
 the corresponding projective scheme is the stable
image of $X$ under the birational map to projective space over $Y$
 given by the linear series 
$|-mD|$ for $m \gg 0$.

 Another viewpoint involves the concept of the {\it
volume} of a  rank one valuation  centered on a local
domain.  By definition, the volume of $\nu$ on $R$ is
$$ {\text{vol}}_R(\nu) := \limsup_{m\ra \infty} \ 
\frac{\text{length} \left(R/\fra_m\right)}{m^n / n!},
$$ where $n$  is the dimension of the local ring $R$.
In  case  $\fra_m$ is  the $m$-th power
 of  a fixed ideal $\fra$,
 then of course this is simply  the multiplicity
of
$\fra$.   But in general the volume is not actually
a multiplicity: indeed, it can be an irrational number; see 
Example \ref{vol.int} (iii).
However Theorem A implies that  the volume of an
Abhyankar valuation is approximated arbitrarily well by
the multiplicities of the ideals
$\fra_m$. Specifically:
\begin{ccorollary.a} If $\nu$ is an Abhyankar valuation
as above, then
\begin{equation*} {\text{vol}}_R(\nu) \ = \ \lim_{m\ra
\infty} \
\frac{e(\fra_m)}{m^n}. \notag
\end{equation*}
\end{ccorollary.a} 
\noi This corollary in turn leads to an interesting
upper bound in the spirit of Teissier
  on 
 the  multiplicity of an
$\frmm$-primary  ideal $I
\subseteq R$   in terms of the  volumes of
its Rees valuations $\nu_1, \ldots,
\nu_r$.\footnote{Recall that the Rees valuations of $I$
are defined by taking   the normalized blowup $X$ of $I$,
and writing 
 $I \OO_X = \OO_X(-e_1E_1 - \dots - e_r
E_r)$ for some 
 prime Weil-divisors $E_i$ of $X$ and some positive
integers
$e_i$.
 Then $\nu_i$ is the valuation on the fraction field of
$R$ given by  order of vanishing along $E_i$.}
Specifically, we  prove an inequality
of Minkowski-type on the volumes of Abhyankar valuations
which implies  that  the multiplicity of $I$ satisfies 
\begin{equation*}
 e(I)^{1/n} \ \leq \ e_1 {\text{vol}}_R(\nu_1)^{1/n} +
\dots +  e_r  {\text{vol}}_R(\nu_r)^{1/n},
\end{equation*} where   
$e_i = \nu_i(I)$,  and as above $n = \dim
R$.
Since an earlier version of this manuscript was written, 
Mircea Mustata has shown that 
 Corollary C and the related Minkowski-type statements
hold for arbitrary valuations (not just Abhyankar) and
even much more generally. See \cite{Must}.
\bigskip

The proof of Theorem A uses the theory of multiplier
ideals.  One can associate to the  valuation ideals $\{ \fra_m \}$
a sequence of \textit{asymptotic multiplier ideals} $\{ \frj_m \} _{m\in
\Phi}$, as defined (in a slightly different setting) in
\cite{ELS}.  The general theory --- which we review in \S 1 --- shows
that these satisfy
\begin{equation} \label{First.Eqn.Thm.Intro}
\fra_m^{\ell} \ \subseteq \ \fra_{m\ell} \ \subseteq \ 
\frj_m^{\ell}
\end{equation}
 for all  
$m \in \Phi$  and  
$\ell \in \NN$.
The main work of the present paper --- which we carry out in \S 2 --- is
to establish that
$\fra_m$ and
$\frj_m$ have ``bounded difference" in the sense that there is a fixed
$\frm$-primary ideal $\frd \subseteq R$ such that 
\begin{equation*} \label{Second.Eqn.Thm.Intro}
 \frd \  \subseteq \ \big( \fra_m : \frj_m \big)
\end{equation*}
for all sufficiently large $m \in \Phi$. Theorem A then follows from 
(\ref{First.Eqn.Thm.Intro}) upon taking $e = \nu(\frd)$. 
 In \S 3, we   discuss the volume of
a valuation: we believe that this is an invariant of
independent interest. Finally we give in \S 4 some
further applications, including the proof of Corollary B.

The present paper continues the
project started in \cite{ELS} of using ideas from
higher dimensional complex geometry to search for
possibly  unexpected uniform behavior in Noetherian
rings. Our techniques rest on  resolution of
singularities and vanishing theorems, and so are
essentially limited to local rings coming from smooth
complex varieties. 
 Although one could expect
the statements themselves to  remain valid in  less
restrictive settings, we haven't seriously investigated
the extent to which such a generalization is possible.  
We hope however that the results appearing here will
pique the interest of experts and encourage them to put
the picture in a broader perspective. 
 We note that the main theorem here was
inspired by an attempt to ``localize" a result of
Fujita concerning the volumes of big line bundles: a
proof via multiplier ideals appears in  
\cite{DEL} (see also \cite{MIAG} and \cite[14.5]{Dem}).
The reader may consult  \cite{Dem}, \cite{Siu},
\cite{nadel}, \cite{Kawamata}, \cite{Demailly.Bourbaki},
\cite{MIAG}, \cite{PAG}, \cite{Ein}  and the 
references therein for numerous  other recent
applications of multiplier ideals to the global study
linear series on a projective variety.

The authors are grateful to Dale Cutkosky, Will Traves and
 especially  Mark Spivakovsky for numerous helpful
and encouraging  discussions. In particular, Cutkosky
first  explained to us a special case of Proposition
\ref{toric}, and Example \ref{vol0}  was worked out in a
series of discussions with   Spivakovsky and 
Traves. Also, we thank Ray Heitmann for questions and remarks
that improved the presentation of our results.

\section{The Asymptotic  Multiplier Ideal of a graded
family}

In this section, we recall the notion of a graded family
of ideals and examine the elementary properties of 
graded families of valuation ideals.  We also  review
the construction and basic properties of the asymptotic
multiplier 
 ideals introduced in  \cite{ELS}.  Because 
 collections of valuation ideals are naturally indexed
by semigroups slightly more general than the  natural
numbers,  it is convenient to  allow graded families 
that are indexed by additive subsemigroups of the real numbers.
 Thus  the exposition here is  slightly more general
than what is  stated in \cite{ELS}, but all the proofs
are the same.

\subsection{Graded Families indexed by a semi-group} Let
$\Phi$ be an additive  subsemigroup of the non-negative real
numbers.
 Of course, one of the main examples    is the semigroup
of natural numbers $\NN$.

\begin{definition} Fix a ring $R$.
 A \textit{graded family} or \textit{graded system of
ideals  indexed by } 
$\Phi$ 
 is a collection $ \fra_{\bullet} = \{ \fra_m \}_{m\in
\Phi}$ of ideals of $R$ satisfying
\begin{equation} \label{gsi.cond} 
\fra_m \cdot \fra_\ell \ \subseteq \ \fra_{m + \ell}
\fall  \ m, \ell \in \Phi. 
\end{equation} To avoid trivialities, we assume also
that $\fra_m \neq 0$ for
$m \gg 0$.
\end{definition} It is convenient to assume also that 
$\fra_0 = R$, in which case condition
 (\ref{gsi.cond}) is  equivalent to the statement that
 the $R$-module
 $$\bigoplus_{m \in \Phi} \ \fra_m
$$ has the natural structure of a $\Phi$-graded
$R$-algebra. We refer to this ring as the {\it  Rees
algebra of the graded system} 
 $\fra_{\bullet}$.
 The theory of asymptotic multiplier ideals is
particularly useful when the Rees algebra fails to be
finitely generated (or at least is not known to be so). 
Of course, one can also define graded families, and
develop the theory of asymptotic multiplier ideals,  for 
 families of coherent ideal sheaves in the structure
sheaf of
 a scheme.  But since all the definitions are local in
nature this involves no essential differences from the
affine setting.

\begin{example}\label{first.examples} The following are
familiar examples of graded families  indexed by the
natural numbers.
\begin{itemize}
\item[(i).] The simplest example is the collection
$\{\fra^m\}$ of powers of a fixed ideal $\fra$. This
should be considered a trivial example of a graded
family. 
\item[(ii).] A slightly less trivial example is the
collection 
$\{\overline{\fra^m}\}$  of integral closures of powers
of  a fixed ideal $\fra$. From the point of view of
multiplier ideals, however, this example is no less
trivial than the first, since the multiplier
ideals are the same. See also Example
 \ref{vol.int} (i).
\item[(iii).] The graded family $\{\fra^{(m)}\}$ of
symbolic powers of a fixed ideal $\fra$  was the main
example treated in \cite{ELS}.
\item[(iv).] The collection $\{\frb_m\}$  of defining
ideals for the base loci of the complete linear series
$|mD|$, where $D$ is a fixed big  line bundle on a
projective variety $X$, forms a graded family of ideals
in the sheaf of rings $\OO_X$. This graded family plays
a central role in Kawamata's work  on  deformation of
canonical singularities \cite {Kawamata}.  
\end{itemize}
\end{example}

\subsection{Graded systems arising from a valuation}

In this paper, we are primarily interested in  graded
families of valuation ideals for  rank one valuations
centered on a regular local domain.  Let us recall some
of the basic terminology and give a few examples. Good
general references on valuation theory are \cite{Vaq}
and \cite{ZS}.

Let $\Gamma$ be an ordered Abelian group, written
additively.
 A valuation on a field $K$ with values in $\Gamma$  is
a  map of Abelian groups
$$
\nu: K^* = K \setminus 0  \lra \Gamma
$$ satisfying the additional condition that
\begin{equation}\label{val.def}
\nu(f + g) \geq \min\{\nu(f), \nu(g)\}.
\end{equation} Since only the image of $\nu$ is of
importance, we assume that $\nu$ is surjective.  An
elementary but useful observation is that in fact,
\begin{equation}\label{val.eq}
\nu(f + g)  = \min\{\nu(f), \nu(g)\} {\text{ whenever }}
\nu(f) \neq \nu(g).
\end{equation}

The set of all elements of $K$ which have non-negative
values  forms of subring $R_{\nu}$, called the {\it
valuation ring} of $\nu$. The valuation ring of $\nu$ is
a local ring with maximal ideal consisting of all the
positively valued elements of $K$, but it is not
Noetherian in general.  One can also define a valuation
ring abstractly as a  subring of $K$ containing either
$x$ or $\frac{1}{x}$ for every
$x \in K$. 
 The data of a valuation ring inside a field
$K$ is equivalent to the data of a valuation  on $K$ (up
order preserving isomorphism of the value group). The
{\it rank}  of a valuation is by definition the Krull
dimension of the valuation ring $R_{\nu}$.

In this paper, we consider only {\it rank one}
valuations on function fields. In this case, the ordered
 group $\Gamma$ can and will be identified with an ordered
subgroup of the real numbers,  the field 
$K$ is assumed to be
 a finitely generated  field extension of some fixed 
ground field $k$, and we consider only those valuations
vanishing on  $k$.  In particular, the valuation ring
$R_{\nu}$ is a $k$-algebra, though not usually finitely generated.

Let $\nu$  be a valuation on a function field $K/k$. The
valuative criterion for properness (see \cite[p.
101]{Hart}) ensures that for  any complete
algebraic variety
$X$ over $k$ with function field $K$ (that is, for any
{\it complete model of }  
$X$),  there is a unique map $\spec R_{\nu} \lra X$.
 Thus  a valuation  chooses,
 in a consistent way, a (not necessarily closed) point
on  every complete model of $K$, namely,  the  the
image  $W$  of the  closed point of $\spec R_{\nu}$
under this map.  The variety{\footnote{Here and
elsewhere, we abuse terminology by failing to
distinguish between an irreducible variety and its 
generic point.}}
  $W$
 is called the {\it center of $\nu$ on $X$}, and the
local ring $\OO_{W,X}$ is called the {\it local ring of} 
$\nu$ on
$X$. The defining ideal $\II_W$ of $W$ consists of all
local sections of $\OO_X$  of positive value.   

For a local domain $R$ contained in $K$, we say that 
the valuation $\nu $  is {\it centered } on $R$  if
$\nu$  takes non-negative values on $R$ and strictly
positive values on the maximal ideal of $R$. Thus the
valuation is centered on any of its local rings.  

Let $\nu$ be a valuation on a function field $K/k$  and
let $X$ be any irreducible $k$-scheme with function
field $K$. For each non-negative $m \in \RR$, the subset 
$$
\fra_m = \{ f \in \OO_X  \, |  \, \nu(f) \geq m\} 
$$ forms an ideal sheaf in  $\OO_X$, called a $\nu$-{\it
valuation ideal} (or just a {\it valuation ideal} when
$\nu$  is understood).  
 When it is necessary to emphasize the model, we will
write $\fra_m(X)$.  Note that  if 
$\pi: X \lra Y$ is a proper birational map between
complete models of $K$,
 then $\pi_* \fra_m(X) = \fra_m(Y).$

 For any subsemigroup $\Psi$ of $\RR$ (eg,
$\Psi = \NN$ or $\Psi = \RR_{\geq 0}$),   the collection
$\{\fra_m\}_{m \in \Psi}$ forms a  graded family of 
ideals in $\OO_X$,
 called the
 {\it graded family of $\nu $ on $X$.\/} It is natural
to index this graded system by $\Gamma$, but allowing
the indexing set to be $\NN$ or $\RR$ will sometimes by
more convenient. The {\it value semigroup}
  $\Phi = \nu(\OO_X) $ 
 is the 'optimal' indexing set: every $\nu$-valuation
ideal  appears exactly once  as a member this graded
family.
 When $X$ is Noetherian, the  value semigroup 
$\Phi$ is well-ordered, meaning that every subset has a
minimal element.

\begin{proposition}\label{primary}
 If $\nu$ is a rank one valuation on a function  field
$K/k$ and $X$  is any complete model of $K$, then the
valuation ideals
$\fra_m(X)$ are primary to the ideal  defining the
center of $\nu$  on
$ X$. If $\nu$  is a rank one valuation centered on
local domain $R$,   then each of the valuation  ideals
$\fra_m$ is primary to the  maximal ideal of $R$.
\end{proposition}

\begin{proof} Suppose that  $f$ and $g$ are local
sections of
 $\mathcal  O_X$ with
$f \notin \fra_m$ but $fg \in \fra_m$. Then 
$\nu(fg) = \nu(f) + \nu(g) \geq m$, but $\nu(f)$ is
strictly less than $m$. This means that $\nu(g)$ is
positive, so some
 positive multiple, say $k\nu(g)$, exceeds $m$. Thus
$g^k $ is in $\fra_m$,  proving that $\fra_m$ is
primary. Furthermore,
 since any element of $\II_W$ has a power in $\fra_m$, 
the radical of $\fra_m$ is  $\II_W$. The argument for
the second statement in the proposition is the same. 
\end{proof}

\begin{example}\label{val.ex} We give four simple
examples of graded families of valuation ideals in the
local ring $R = k[x,y]_{(x, y)} $  of the origin in the
affine plane. Each arises from a different rank one
valuation on the function field $K = k(x, y)$.
\begin{itemize}
\item[(i).] Let $\nu$  be the valuation given by 'order
of vanishing at the origin'. Explicitly, for a
polynomial $f$,  $\nu(f)$ is the degree of the smallest
degree  non-zero monomial appearing in the unique
expression of $f$ as a sum of monomials $x^ay^b$. The
value of any rational function $\frac{f}{g}$ where $f$
and $g$  are  polynomials is uniquely determined by
virtue of  $\nu$  being a group homomorphism:
$\nu(\frac{f}{g}) = \nu(f) - \nu(g)$.  

This valuation has value group $\ZZ$ and value semigroup
$\NN$ on $R$.
 The valuation ideals of $\nu$  are powers of the
defining ideal of the origin:
$$
\fra_m = (x, y)^m.
$$
\item[(ii).] Let $\pi: X \lra \spec{R}$  be any proper
birational map from a normal scheme $X$  and let $D$  be
any prime divisor of $X$  collapsed to the origin. Let
$\nu_D$  be the valuation given by  'order of vanishing
along $D$.' Explicitly, for any $f \in R$,  consider the
image of $\pi^*f$ in the local ring $\OO_{D,X}$ of
$X$ at the generic point of $D$. Then $\nu_D(f)$ is the
maximal integer $n$  such that $\pi^*(f)$ is divisible
by $t^n$ in $\OO_{D,X}$, where $t$ is a uniformizing
parameter for the discrete valuation ring 
 $\OO_{D,X}$.

This valuation has value group $\ZZ$, and value
semigroup on $R$  a subsemigroup of  $\NN$. The
valuation ideals are 
$$
\fra_m = \pi_* \OO_X(-mD),
$$ which in general, can be tricky to understand. In the
simple case where $\pi$ is the blowup of the origin and
$D$ is the resulting exceptional divisor, this example
recovers Example (i) above. 
\item[(iii).]  Let $\nu$  be the valuation on $k(x, y)$ 
 defined by  the assignment $\nu(x) = 1$  and $\nu(y) =
\pi \in \RR$. This uniquely determines a valuation on
$k(x, y)$: 
 the value of each monomial $x^ay^b$ is $a + b \pi$, and
because $1$ and $\pi$ are
 $\ZZ$-independent, distinct monomials have distinct
values, so
 the value of an arbitrary  polynomial is determined
 by  (\ref{val.eq}).

This valuation has value group $\ZZ + \ZZ \pi \subset
\RR$ and value semigroup  $\NN  + \NN \pi$ on $R$.  The
valuation ideals for $\nu_{\pi}$ are all monomial ideals
$$
\fra_m = (\{x^ay^b  \, | \,  a + \pi b \geq m\} ). 
$$
\item[(iv).] Let $\nu$ be the valuation given by 'order
of vanishing along  the analytic arc $y = e^x - 1$'.
Explicitly, for a polynomial $f$,  we define $\nu(f)$ to 
be the highest power of $t$ dividing the image of $f$
under  the map
$$ k[x, y] \hookrightarrow k[[t]]
$$
$$ f(x, y) \mapsto f(t, e^t-1).
$$ Here $e^t-1$ denotes the power series $t +
\frac{t^2}{2!} + \frac{t^3}{3!} + \dots$ (assuming that
$k$ has characteristic zero). 

This valuation has value group $\ZZ$ and  value
semigroup  $\NN $ on $R$. Its valuation ideals are given
by
$$
\fra_m = (x^m, y - x - \frac{x^2}{2!} - \dots -
\frac{x^{m-1}}{(m-1)!}).
$$
\end{itemize}
\end{example}

The preceding examples are  simple examples of rank one
valuations centered on the origin of  the plane.  The
first three are {\it Abhyankar valuations} (see
\ref{ZA}); all four have finitely generated value
groups. By contrast, there are valuations with value
group $\QQ$  centered on the origin of the plane (see
Example \ref{vol0}); while more complicated to describe,
 these are actually 'typical' in a   certain sense.
 The classification of valuations centered on the origin
of the plane is a beautiful story; see \cite{Spiv}.  
For further examples of valuations, including higher
rank valuations, consult \cite{Vaq} \S 10,  or
\cite{ZS} VI \S 15.
\endexample

\subsection{Multiplier ideals.}

We now recall the construction and basic properties of
multiplier ideals as well as the asymptotic
constructions from 
\cite{ELS}. We will give only a few proofs; the rest can be found
for instance in \cite{ELS}, \cite{MIAG} or \cite{PAG}.

Consider a scheme $X$, smooth and  essentially of
finite type over a field of characteristic zero---
mainly we have in mind the case where $X$ is a smooth complex variety or
the spectrum of a regular local $k$-algebra.  Let
$\fra
\subseteq \OO_X$ be a (coherent) ideal sheaf on $X$.
Given a rational number $c > 0$ we  define
\textit{multiplier ideals} 
\[  \MI{X; c \cdot \fra} \ = \ \MI{X; \fra^c} \
\subseteq \ \OO_X. \]
Intuitively these are ideals with remarkable
cohomological properties which reflect in a somewhat
subtle manner  the singularities of the divisors of
functions $f \in \fra$. The construction starts by
taking  \textit{log resolution}
$\mu: X^\pr
\lra X$
 of $\fra$. Recall that this means that $\mu$ 
is 
 a projective birational map from a regular scheme
$X^\pr$  to $X$
 such that $\fra \OO_{X^\pr} =
\OO_{X^\pr}(-F)$, where   $F$ is an effective Cartier
divisor on
$X^\pr$ with the property that the sum of $F$ and the
exceptional divisor of $\mu$ has simple normal crossing
support. Such resolutions can be constructed (as we are
in characteristic zero)  by resolving the singularities
of the blow-up of $\fra$. We write $K_{X^\pr/X} =
K_{X^\pr} -
\mu^* K_X$ for the relative canonical divisor of
$X^\pr$ over $X$. Given $\fra$ and $c >0$ as above, we
now define
\begin{equation}\label{mult.def}
\MI{X ; c \cdot \fra} \ = \ \MI{c \cdot \fra} \ = \
\mu_*
\OO_{X^\pr}\big(\, K_{X^\pr / X}  - [ \, c
\, F\, ]\, \big ).
\end{equation}
 Here $c F$ is viewed as an effective
$\QQ$-divisor on $X^\pr$, and its integer part $[ c \,F]$
is defined by  replacing the coefficient of each
component by the greatest integer less than or equal to
it.  This definition is independent of the log
resolution $\mu$; see {\it e.g.} \cite{MIAG}. When $X$
is the affine  scheme $\spec R$, we write also 
$ \MI{R , c \cdot \fra},$ and when $c = 1$ we write
simply $\MI{\fra}$.

\begin{remark} Following Lipman \cite{Lipman} it is also
possible to define
$\MI{X ; c
\cdot \fra}$ without referring to a log resolution. For
example, when $R$ is a regular local ring,
\begin{equation}
\MI{ R ; c \cdot \fra} = \bigcap_{\nu} \{r \in K \, | \,
\nu(r) 
\geq c \nu(\fra) - \nu(Jac_{R_{\nu}/R}) \},
\end{equation} where the intersection is taken over all
valuations of $K$ given by order of vanishing along some
prime divisor on a model $X^\pr$  of $K$, and where 
$Jac_{R_{\nu}/R}$ denotes the Jacobian ideal of the
extension 
$R \hookrightarrow R_{\nu}$. This is a slight
modification of Lipman's definition of an adjoint ideal,
in which  we have allowed for the possibility of a
coefficient $c$; see \cite{Lipman}.  This approach
makes sense for arbitrary  regular Noetherian  schemes
(not just for classes of schemes that  admit a good
theory of resolution of singularities as we have assumed
here.) However  multiplier ideals derive their power
from the properties they satisfy, and the important
facts --- which rest on vanishing theorems --- are so
far only known over fields of characteristic zero. 
\end{remark}

\begin{remark} Multiplier ideals were originally defined
analytically, as    ideals of germs of holomorphic
functions that are $L^2$-integrable for a certain
weighted $L^2$ space. See  {\it e.g.} \cite{Dem}. In
this approach they appear as sheaves of multipliers,
whence the name. 
\end{remark}

Among the various properties these ideals satisfy, we
will need three  in particular. First:
\begin{equation} \label{Prop.MI.1}
\fra \ \subseteq \ \MI{\fra} \end{equation}
for any ideal $\fra$. This is elementary: it boils down
to the fact that the relative canonical bundle
$K_{X^\pr/X}$ of a log resolution is effective.
More substantially for any rational $c > 0$ and any
$\ell \in \NN$ one has the
\textit{subadditivity  relation}:
\begin{equation} \label{Prop.MI.2}
\MI{ \, \ell c \cdot \fra \,} \ \subseteq \
\MI{\,
 c \cdot \fra 
\,}^\ell.
\end{equation}
This is established in \cite{DEL} using vanishing
theorems. The third  useful property, which follows easily
from the definitions, is the behavior of these ideals
under birational maps. Specifically, if $\mu : X^\pr
\lra X$ is a proper birational map, with $X^\pr$ regular,
and if $\fra \subseteq \OO_X$ is any ideal, then
\begin{equation} \label{birat.transf.rule}
\MI{X ; \fra} \ = \ \mu_* \Big( \MI{X^\pr ; \fra^\pr}
\otimes \OO_{X^\pr}(K_{X^\pr / X}) \Big),
\end{equation}
 $\fra^\pr = \fra \cdot \OO_{X^\pr}$ being the
pullback of $\fra$ to $X^\pr$. The reader may consult
\cite{MIAG}, or \cite[Part III]{PAG} for proofs of these
and other properties of multiplier ideals. 

Given now a graded family $\fra_{\bullet}$ as above and
an index $m \in \Phi$, we will construct an
\textit{asymptotic multiplier ideal} $\frj_m =
\frj_m(\fra_{\bullet})$ which reflects the asymptotic
properties of all the ideals $\fra_{p m}$ for $p \in
\NN$. From  (\ref{gsi.cond}) and (\ref{mult.def}), it is
easy to check that for each
$m \in \Phi$,  we have  
\begin{equation} \label{gri.inclusions}
\MI{\fra_{m}} \subseteq
\MI{\tfrac{1}{p} \cdot \fra_{pm}} 
\end{equation} for all $p \in \NN$; see \cite[\S1]{ELS}.
This, together with the Noetherian property for $\mathcal
O_X$, implies that  the set of ideals 
\begin{equation}  \Big \{
\MI{\tfrac{1}{p}
\cdot \fra_{pm}} \Big \}_{p  \in \NN} 
 \notag
\end{equation} has a unique maximal element. We then
define  the $m$-th  asymptotic multiplier ideal
$\frj_m(\fra_{\bullet})$ to be this maximal element. 
 In other words,  
\begin{equation}  \label{def.ami.eq} 
\frj_m(\fra_{\bullet}) 
\ = \ \MI{\tfrac{1}{p} \cdot   \fra_{pm}}
\ \for \text{sufficiently divisible } p \in \NN.
\end{equation} In fact, it is not necessary to assume
that 
$p$ is  sufficiently divisible:  assuming that $\fra_m
\ne (0)$ for $m \gg 0$,  any very large
$p$ will do;  
 see \cite{ELS}, Remark following Definition 1.4. 

The essential property of these ideals is summarized in
the next result, which was established in
\cite{ELS}.
\begin{theorem} \label{Properties.Mult.Ideal.Thm}
 For any graded system $\fra_{\bullet}$,
 any index $m \in \Phi$, and any natural number $\ell
\in \NN$ one has inclusions:
\begin{equation} \label{Mult.Ideal.Thm.Eqn}
\fra_{m}^\ell \ \subseteq \ \fra_{m \ell} \ \subseteq \
\ \frj_m^\ell. 
\end{equation}
\end{theorem}
\begin{proof}[Sketch of proof] The first inclusion is
definitional. For the second, note first that $\fra_{m
\ell} \subseteq \frj_{m \ell}$: this is easily checked
using (\ref{Prop.MI.1}) and (\ref{gri.inclusions}). So
it is enough to show that $\frj_{m \ell} \subseteq
\frj_{m}^\ell$. For this, choose any $p \gg 0$. Then
using the subadditivity relation (\ref{Prop.MI.2}) one
finds:
\begin{align*}
\frj_{m \ell} \ &= \ \MI{\tfrac{1}{p} \cdot \fra_{m \ell
p}}\\
			&= \ \MI{ \tfrac{\ell}{\ell p} \cdot \fra_{m \ell
p}} \\
&\subseteq \  \MI{ \tfrac{1}{\ell p}\cdot \fra_{m \ell 
p}}^\ell
\\
&= \ \frj_m^\ell, 
\end{align*}
as required. \end{proof}

\begin{remark} The second inclusion in
(\ref{Mult.Ideal.Thm.Eqn}) does not hold in general if
one works with the ``absolute" multiplier ideal
$\MI{\fra_m}$ in place of $\frj_m$. 
\end{remark}

\section{Abhyankar Valuations}

Let $\nu$ be a rank one valuation on a function field
$K/k$. There are two basic invariants of $\nu$. The {\it
rational rank }  of $\nu$ is the dimension of the
$\QQ$-vector space $\QQ \otimes_{\ZZ} \Gamma$.  The
{\it  transcendence degree}  of $\nu$ is the
transcendence degree of the residue field of the
valuation ring $R_{\nu}$  over
$k$. Equivalently, the transcendence degree is
 the maximal dimension of the center of $\nu$  over all 
models of $K/k$. The basic result relating these
invariants is the Zariski-Abhyankar
inequality.{\footnote{In the form stated, Theorem
\ref{ZA} is due to Zariski;
 see \cite{ZS}.  Abhyankar \cite{Abh} later proved a
more general version of \ref{ZA} for valuations centered
on any local Noetherian  domain (the non-function  field
case).}} 

\begin{theorem}[The Zariski-Abhyankar
Inequality]\label{ZA} For any valuation on a function
field $K/k$
\begin{equation}
\label{ab} {\textnormal{trans.deg }} \nu +
{\textnormal{rat.rank }} \nu \leq \dim K/k.
\end{equation} Furthermore, if equality holds in
(\ref{ab}), then the value group
 $\Gamma$ is a finitely generated (free) Abelian group.
Here, $\dim K/k$ refers to the transcendence degree of
$K$ over $k$, or equivalently, to the  dimension of any
complete model for $K/k$.
\end{theorem}

A valuation satisfying equality in (\ref{ab}) is called
an
 {\it Abhyankar valuation.} Abhyankar valuations
generalize the familiar example of  divisorial
valuations, that is, valuations given by order of
vanishing along some divisor on a normal model of $K/k$. 
Note that divisorial valuations have rational rank one
and transcendence degree $n-1$, where $n$ is the
dimension $\dim K/k$. 

\begin{example} For the valuations in  Example
\ref{val.ex}, the rational rank and transcendence degree
are easily computed:  for (i) and (ii),
 the rational rank is 1 and the transcendence degree is
1; for (iii), the rational rank is 2 and the
transcendence degree is zero; and for (iv), the rational
rank is 1 and the  transcendence degree is zero.
\end{example}

We  now state the main technical result of the present
paper:

\begin{theorem}\label{main} Let $\nu$ be an Abhyankar
valuation
 on a function field $K/k$   of characteristic zero. Let
$R$ be the local ring of the center of $\nu$ on some
smooth  model of $K$, and let $\{\fra_m\}_{m\in \Phi}$
 be the graded family of $\nu$-valuation ideals in $R$.  
Then there exists a non-zero element $\delta \in R$ such
that 
\[ \delta \cdot \frj_m \ \subseteq \  \fra_m\] for all $m
\in
\Phi$, where $\{\frj_m\}_{m\in \Phi}$
 are the  asymptotic multiplier
ideals associated to $\fra_{\bullet}$.  
\end{theorem}

\begin{remark} For any $\nu$-valuation ideal $\fra$ in
domain $R$ and for an arbitrary ideal $\frb$,  the colon
ideal $(\fra: \frb)$ is also a $\nu$-valuation ideal of
$R$; see
\cite[p. 342]{ZS}. Thus (taking $\Phi$ to be the full
value semigroup $\nu(R)$), 
 there is a function $\alpha: \Phi \lra \Phi$ such that 
$$ (\fra_m: \frj_m) = \fra_{\alpha(m)}.
$$ Theorem \ref{main} says that for Abhyankar valuations
centered on a regular local domain essentially of finite
type over a field of characteristic zero, this function
is bounded above.  In other words, the graded family of
($\frak m$-primary) valuation ideals $ (\fra_m: \frj_m)$ has a minimal
element. 
\end{remark}

In view of the preceding Remark,  Theorem A 
from the Introduction follows immediately
upon combining  Theorems
\ref{Properties.Mult.Ideal.Thm} and \ref{main}.  
Indeed,  one simply takes
 $e = \nu(\delta)$,
where $\delta$ is the non-zero element of $R$ whose existence is
guaranteed by Theorem \ref{main}. Since $\delta \frak j_m \subset \frak a_m$
for all $m$,
clearly $\frak j_m \subset \frak a_{m - e}$ for all $m$,
and so $\fra_m^{\ell} \subset \frak j_m^{\ell} \subset
 \fra_{m-e}^{\ell}.$

Theorem \ref{main} also implies that Izumi's theorem
holds for non-divisorial Abhyankar valuations in our
setting:
\begin{corollary}[Izumi's Theorem for Abhyankar valuations] \label{izumi} Let
$\nu$ and $w$  be rank one Abhyankar valuations on a function
field $K/k$  of characteristic zero and let $(R, \frmm)$
be any regular local ring essentially of finite type
over $k$ on which  both $\nu$ and $w$ are centered. 
Then there exists 
  $C>0$ such that
\begin{equation}\label{iz}
\nu(x) \geq C w(x)
\end{equation} for all non-zero elements $x \in R$.
\end{corollary}

\begin{proof}
Without loss of generality, we may assume that 
  the minimal value  obtained by $w$ 
on $R$ is 1.  Indeed, 
because $R$ is Noetherian, $w$ must achieve some minimal
value on $R \setminus 0$;  now we can simply scale the  
values of $w$ and $\nu$ by this minimal value.

Now we claim that 
there 
exists a value $p$ such 
that 
\begin{equation}\label{linear}
\fra_{mp} \subset \frb_{m}
\end{equation}
for all $m \in \NN$, where $\{\fra_{m}\}$ (respectively $\{\frb_m\}$)
denote  valuation ideals of $\nu$ (respectively, $w$).
Indeed, note that 
$\frb_1 = \frm$, and so 
$\frmm^{\ell} \subset \frb_{\ell}$
for all $\ell \in \NN$. 
Thus to prove (\ref{linear}), 
it is enough to show
that 
there exists $p$ such that 
 \begin{equation}\label{max}
\fra_{p\ell} \subset \frmm^{\ell} 
\end{equation}
 for all $\ell$.
(In other words, we are reduced to the case where $w$ is the $\frmm$-adic
 valuation on $R$.)
By Theorem
\ref{Properties.Mult.Ideal.Thm}, we see that 
(\ref{max}) follows  immediately provided that 
 some $\frj_p \subset \frmm$. But if 
  all $\frj_p$ are trivial, we have  $\delta \in  
\cap \fra_m = (0)$, contradicting Theorem \ref{main}.

Finally, (\ref{iz}) follows
as in \cite[Lemma 1.4]{Hub.Swan}
 by setting 
  $C = 2p -1$ (and enlarging $p$ if necessary so that $p \geq 2$).
Indeed, 
suppose on the contrary that  there is some $x \in R$ such that 
$\nu(x) > Cw(x)$. Set $m = w(x) + 1.$
Then $x \in \fra_{mp}$, but  $x \notin \frb_m $.

\end{proof}

\begin{remark}\label{fail} Theorems A and \ref{main}
--- and also the Izumi-type statement of Corollary \ref{izumi} ---
can fail for non-Abhyankar valuations. In fact, consider
the  valuation $\nu$  given by order of
vanishing along the exponential curve $y = e^x - 1$  
(Example \ref{val.ex}.iv). Here $\fra_{p} \subseteq R
= \CC[x,y]_{(x,y)}$ has colength $p$. So there cannot
exist a non-trivial ideal $\frj \subseteq R$ having the
property that $\fra_{m \ell} \subseteq \frj^{\ell}$ for
fixed $m$ and $\ell \gg 0$, since the colength of
$\frj^\ell$ would grow quadratically in $\ell$.
Therefore the inclusion (\ref{First.Eqn.Thm.Intro}) in
Theorem A can only hold with 
$\frj_m = R$ for all $m$. On the other hand,
$\cap_m \fra_m = (0)$, so there cannot exist a fixed
non-zero element $\delta \in R$ with $\delta \cdot
\frj_m \subseteq \fra_m$ for all $m$. However for an
arbitrary valuation $\nu$, it is possible that the
colon ideals $\frd_m = (\fra_m : \frj_m)$ ``grow
slowly": see Remark
\ref{tightbound} for a precise
statement. 
\end{remark}

We now prove Theorem \ref{main}. The outline is this:
Lemma \ref{blowup} below guarantees that   it is enough
to  find  $\delta$ after blowing up,   but after  a
suitable blow up, any Abhyankar valuation is
``essentially monomial''  by Proposition \ref{toric},
where a direct computation can be carried out.

\begin{lemma}\label{blowup} Let $\pi: X \lra Y$ be  a
proper birational map of smooth varieties over $k$.
Assume that there exists an ideal $\frd^\pr  \subseteq
\mathcal  O_X$ such that
$\frd^\pr \cdot \frj_m(X) \subset \fra_m(X)$.
    Then  there exists an ideal $\frd \subseteq
\mathcal  O_Y$ such that $\frd \cdot   \frj_m(Y)
\subseteq
\fra_m(Y)$.
\end{lemma}

\medskip
\begin{proof} Let $\frd = \pi_* \big(\frd^\prime
\omega_{X/Y}^{-1}\big)$.  Note that both $\frd^\prime$
and
$\omega_{X/Y}^{-1}$ are ideals of $\mathcal  O_X$,  so
$\frd$ is an ideal of $\mathcal  O_Y$.
Now
$$
\frd \cdot  \frj_m(Y)\  = \  \pi_*\big(\frd^\prime
\omega_{X/Y}^{-1}\big)
\cdot \pi_*\Big(\MI{\, X,
\frac{1}{p}\cdot \fra_{mp}(Y)\cdot \mathcal O_X \,}
\otimes
\omega_{X/Y}\Big)$$ where $p$ is sufficiently large.
Here we have used the definition of the asymptotic
multiplier ideal together with the transformation rule
(\ref{birat.transf.rule}) for multiplier ideals under
proper birational  morphisms.  
Therefore
\begin{align*}
\frd \cdot  \frj_m(Y) \  &\subseteq \ 
\pi_*\Big(\frd^\pr\omega_{X/Y}^{-1} \cdot \mathcal  J(X,
\frac{1}{p} \cdot \fra_{mp}(Y) \mathcal  O_X) \otimes
\omega_{X/Y}\Big) \\ &= \  \pi_* \Big(\frd^\pr
\cdot \MI{X,
\frac{1}{p}\cdot \fra_{mp}(Y)\cdot \mathcal 
O_X}\Big).
\end{align*}
Because $
\fra_{mp}(Y)\mathcal  O_X \subset \fra_{mp}(X)$ for all
$m$ and $p$, the corresponding inclusion holds for  the
multiplier ideals. 
Note that the multiplier
ideal appearing on the right is that associated to the
pull-back of a valuation ideal on $Y$ rather than the
corresponding valuation ideal on $X$. However
$\fra_{mp}(Y)\cdot \OO_X \subseteq \fra_{mp}(X)$ and
consequently for  $p\gg 0$:
$$\frd^\pr \cdot\MI{X,
\frac{1}{p}\fra_{mp}(Y)\cdot \mathcal  O_X} \ 
\subseteq \ 
\frd^\pr \cdot \MI{X,
\frac{1}{p}\fra_{mp}(X)} \  = \ \frd^\pr \cdot
\frj_m(X).
$$
But $\frd^\pr \cdot \frj_m (X) \subseteq \fra_m(X)$,
and putting these inclusions together we find
that  
\[ \frd \cdot \frj_m(Y) \ \subseteq \ 
\pi_*\big(\frd^\pr
\cdot
\frj_m(X)\big)
\ \subseteq \ \pi_*(\fra_m(X)) \ = \ \fra_m(Y) \] for
all
$m$, as required. 
\end{proof}

The next Proposition is probably well-known, at least
for valuations of transcendence  degree zero.
We learned that case from Dale Cutkosky.
 
\begin{proposition}\label{toric} Let $\nu$ be a rank one
Abhyankar valuation on a function field 
$K/k$ of characteristic  zero.  Given any model $Y$ of
$K/k$, there exists a smooth model $X$ dominating $Y$
and a regular system of parameters 
$x_1, \dots, x_r$  for the local ring of $\nu$ on 
$X$ such that $\nu(x_1), \dots, \nu(x_r)$
 freely generate the value group $\Gamma$.
\end{proposition}

\begin{proof} Let $r$ denote the rational rank of $\nu$,
so  by Theorem \ref{ZA},  
 $\Gamma
\cong \ZZ^r$.  Fix $f_1, \dots, f_r$ in the field $K$,
whose values generate
$\Gamma$.  By replacing $f_i$ by $\frac{1}{f_i}$ if
necessary, 
 we can assume all $v(f_i) > 0.$ 

We can write each $f_i$ as a fraction $\frac {a_i}{b_i}$,
where the $a_i$ and $b_i$ are regular on some
neighborhood of the center of
$\nu$  on  $Y$. By blowing up the ideals $(a_i, b_i)$,  
we can make the fractions $\frac{a_i}{b_i}$ regular on
some neighborhood of
 the center. 
 By blowing up further if necessary, we can assume that
the dimension of the center is the  transcendence degree
of $\nu$--- which means  its codimension equals the 
rational rank in the presence of the Abhyankar
hypothesis.
 So we have created a
 model $Y'$ dominating $Y$ where the elements $f_i$ are
regular on a neighborhood of the center  of $\nu$, and
where the codimension of the center is  exactly $r$, the
rational rank of $\Gamma$.

Now we use embedded resolution of singularities to
resolve the hypersurface defined by the product $f_1
f_2\dots f_r$ in a neighborhood of the center on $Y'$.
This produces for us a smooth model $X$ dominating $Y'$
such that the pullback of the hypersurface to this model
has  simple normal crossing support. In particular, for
any  closed point $x$ of $X$, we have
$$ f_1 f_2\dots f_r = u \, x_1^{a_1} x_2^{a_2} \dots
x_N^{a_N},
$$ where $x_1, \dots, x_N$ is a regular system of
parameters at  $x$,
 the exponents $a_i$ are natural numbers, and $u$ is a
regular function invertible in a neighborhood of $x$. 
Because the local rings of $X$ are unique factorizations 
domains, for each $f_i$ we have
$$ f_i =  u_i \,  x_1^{a_{i1}} x_2^{a_{i2}} \dots
x_N^{a_{iN}}
$$ for some $a_{ij} \in \NN$ and some unit $u_i$.

In particular, choosing the point $x$ to be in the center
 $W$ of $\nu $  on $X$, then the elements $u_i$ are also
units  in the local ring  $\mathcal  O_{W,X}$.
            Because units in $\mathcal  O_{W,X}$ have
value zero, we see
that                                       
$$
\nu(f_i) = \sum_{j=1}^N  a_{ij} \nu(x_j).
$$  So clearly,  the elements
$\nu(x_j)$ generate $\Gamma$.

We claim that exactly $r$ of the elements $x_j$  have
non-zero value. Indeed, if fewer have non-zero value,
then the rank of $\Gamma$ can not be
 $r$. But if more have non-zero value, then there are at
least $r+1$ of the parameters
 $x_1, \dots, x_{r+1}$ contained in the defining ideal
of the  center
$W$.  This would force $W$ to have codimension greater
than $r$, 
 a contradiction.

Relabeling so that the parameters $x_1,  \dots,  x_r$ are
those  with positive value, note finally that the images
of these elements  generate the maximal ideal in the
local ring of $X$  along $W$. Indeed, this maximal ideal
is generated by the  image of the defining ideal $\II_W$
of $W$, and we have already remarked that 
$(x_1, \dots, x_r) \subset \II_W$. But since 
$x_1, \dots, x_r$ are part of regular sequence of
parameters in  a neighborhood of $W$, they must generate
the maximal ideal after localizing at $\II_W$.   Thus
the proposition is proved: the elements $x_1, \dots,
x_r$ of $K$ are a regular system of parameters for the
local ring $\mathcal  O_{W, X}$ and the values 
$\nu(x_1), \dots, \nu(x_r)$ generate $\Gamma$.
\end{proof}

We can now finish the proof of the main theorem. 
\begin{proof}[Proof of Theorem \ref{main}] By Lemma
\ref{blowup}  and Proposition \ref{toric},
 we can assume that we are in the following situation.
The variety $X$ is smooth, and the center $W$ 
 of $\nu$ on $X$ is of
 codimension $r$ equal to the rational rank of $\nu$; 
furthermore, the local ring $R$ of $X$ along $W$ has
regular system of parameters $x_1, \dots, x_r$ whose
values generate the value group $\Gamma$. We wish to
prove that there exists $\delta \subset \mathcal  O_X$
such that
$\delta \frj_m \subset \fra_m$.  Because  $\fra_m$ is
primary to $\II_W$ (see Proposition \ref{primary}),
 it is enough to check  this after localizing along the
defining ideal of $W$,  so we consider the graded family
of valuation ideals $\{\fra_{\bullet}\}$ in  the local
ring $(R, \frak m$.

Because the values of the parameters $x_1, \dots, x_r$
are all $\ZZ$-independent, the ideals $\fra_m$  are 
generated by
 'monomials' in $x_1, \dots, x_r$. Indeed, fix any $m
\in \Phi$.  By (\ref{primary}), some power of the
maximal ideal of $R$, say $\frmm^t$, is contained in
$\fra_m$.   Now consider an 
 arbitrary element $f$ of $R$ (not already in $\frmm^t$).
Modulo $\frmm^t$, $f$ can be  written as a sum  of 
monomials in the regular system of parameters $x_1,
\dots, x_r$ with unit coefficients. Because the  values
of the $x_i$ are independent,
 each of these monomials has a distinct value, and so
the value of $f$ is equal to the value of
 the unique smallest value monomial in this sum (Cf
(\ref{val.eq})).  So each of the 
  monomial ``terms" of $f$ are in  $\fra_m$, and 
$\fra_m$  is generated by monomials in the regular
system of  parameters for $R$.

We now claim that this setup is sufficiently close to the
standard monomial case so as to be able to apply the 
computation derived in \cite{howald} for the multiplier
ideal  of a monomial ideal in a polynomial ring. Roughly
this reason is this: the ring $R$ is \'etale over a
polynomial ring and the computation of multiplier ideals
commutes with \'etale extension. We justify this
carefully in the next paragraph.

Think of the parameters $x_1, \dots, x_r$ as local
sections of $\OO_X$ and extend them to a full set of
regular parameters in 
 some affine neighborhood $U$ of $W$ on $X$. The
inclusion
\begin{equation}\label{incl}
 k[x_1, \dots, x_r, x_{r+1}, \dots, x_n] 
 \hookrightarrow \OO_{X}(U)
\end{equation} induces a natural map 
\begin{equation}\label{map} U \lra \A^n
\end{equation}
 consisting of an open immersion followed by a finite
map.  Localizing at the prime ideal of $W$ and its
corresponding contraction to  the polynomial ring, we
have an inclusion
\begin{equation}\label{incl2}
 A = k[x_1, \dots, x_r, x_{r+1}, \dots, x_n]_{(x_1,
\dots, x_r)} 
 \hookrightarrow \OO_{W,X} = R.
\end{equation} Our claim above that $\fra_m$ is a 
monomial ideal is tantamount  to saying that $\fra_m$ is
the expansion of a monomial ideal $\fra_m' \subset 
k[x_1, \dots, x_n]$. The monomial ideal $\fra_m'$ is
itself a valuation ideal for the valuation on $k(x_1,
\dots, x_n)$ obtained by  restricting
$\nu$ to this subfield. Because the maximal ideal of $A$
expands to the maximal ideal of $R$, the map of rings
(\ref{incl2}) is \'etale, which is to say,  the morphism
(\ref{map}) is \'etale is a neighborhood of $W$
\cite{Milne}. Thus replacing $U$  by a possibly smaller
open neighborhood, we can assume the morphism 
$$ U \lra \A^n
$$ consists of  compositions of open immersions with  a
finite \'etale map.  But the computation of multiplier
ideals commutes with pullback under both open immersions
(obvious) and  finite \'etale maps (straightforward; see
\cite{MIAG}, \S 5.4). So
$$
\MI{ \A^n,  c \cdot \fra_m'} \cdot \OO_U \ = \ \MI{U,  
c
\cdot
\fra_m'\OO_U},
$$ and  passing to the local ring at $W$ (which after
all, amounts to taking a limit of pullbacks to smaller
and smaller affine neighborhoods of $W$), we see that 
$$
\MI{\A^n,  c \cdot \fra_m'} R = \MI{R,   c \cdot \fra_m'
R} = 
\MI{R,   c \cdot \fra_m}.
$$  So  to compute the multiplier ideal of $\fra_m$ in
the local ring $R$, it is sufficient to compute the
multiplier ideal of the monomial ideal
$\fra_m'$ in the polynomial ring $k[x_1, \dots, x_n]$
and expand to $R$.

We recall the formula for the  multiplier ideal of a
monomial ideal from \cite{howald}. Let $\fra  \subset
k[x_1, \dots, x_n]$ be an ideal generated by monomials
and let $L$ denote its lattice of exponents:
$$ L = \{(a_1,  \dots, a_n) \, | \, x_1^{a_1}\dots
x_n^{a_n} \in \fra\}.
$$ Then the multiplier ideal $\MI{c\cdot \fra}$ is  the 
ideal of the polynomial ring  generated  by those
monomials $x_1^{b_1}\dots x_n^{b_n} $ satisfying 
$$ (b_1, \dots, b_n) + (1, \dots, 1) \in
\{{\text{hull}}(cL)\}^{int},
$$  where 
$\{{\text{hull}}(cL)\}^{int}$  denotes the interior of
the convex hull of the  lattice $L$ scaled by the real
number $c$.

In our case, the monomial ideals $\fra_m'$ are generated
by
$$
\{x_1^{a_1}\dots x_n^{a_n} \, | \, \sum a_i \nu(x_i)
\geq m\},
$$ so the multiplier ideals are given by 
$$
\MI{c \cdot \fra_m'} = 
\left(\{x_1^{b_1}\dots x_n^{b_n} \, | \, \sum (b_i +1)
\nu(x_i) >  c m\}\right).
$$ In particular, for any real positive $p$,
$$
\MI{c \cdot \fra_m'} = \MI{\frac{c}{p} \cdot \fra_{mp}'} 
$$ so that the asymptotic multiplier ideals $\frj_m'$ of
the graded family $\{\fra_m'\}$ in the polynomial ring
satisfy  
$$
\frj_m' = \MI{\fra_m'}.
$$ Expanding to $R$, we see that the asymptotic
multiplier ideals in $R$ are given by
$$
\frj_m(\fra_{\bullet}) = 
\left(\{x_1^{b_1}\dots x_n^{b_n} \, | \, \sum (b_i +1)
\nu(x_i) > m\}\right).
$$

Finally, using this description of the  asymptotic
multiplier ideals,  we observe that the element $\delta
= x_1\dots x_r$ in $R$ satisfies the condition that 
$$
\delta \frj_m \subseteq \fra_m
$$ and  Theorem \ref{main} is proved.
\end{proof}

\section{The Volume of a Graded Family of
$\frmm$-primary ideals}

In this section we introduce the volume of a graded
family of $\frakm$-primary ideals, and compare it to the
multiplicities of the individual ideals of the graded  system.
 In particular, we deduce  Corollary
C from the Introduction as a special case of a general phenomenon.

\begin{definition} \label{vol} Let $\fra_{\bullet} =
\{\fra_m\}_{m \in \Phi}$ be a  graded family of
$\frmm$-primary ideals in a  local Noetherian ring 
$(R, \frmm)$ of dimension $n$. The volume of  
  $\fra_{\bullet}$ is the real number
\begin{equation} \label{vol.lim}  \vol(\fra_{\bullet})
\ = \ \limsup_{m \in \Phi} \, 
\frac{\length (R/\fra_m)}{m^n / n!}.
\end{equation}
\end{definition}

\begin{remark} The volume of a graded system is the
local analogue of the volume of a big divisor $D$ on a
projective variety
$X$ of dimension $n$, which is defined to be 
\[ \vol(D) \ = \ \limsup_m
\frac{\hh{0}{X}{\OO_X(mD)}}{m^n/n!}.\]
When $X$ is smooth and $D$ is   ample, this
coincides up to constants with the volume of $X$
determined by any K\"ahler form representing $c_1\big(
\OO_X(D) \big)$, which explains the terminology. 
\end{remark}

While Definition \ref{vol} works well for 
graded systems indexed by the natural numbers $\NN$, it
does not have very good behavior for arbitrary graded
systems $\big \{ \fra_m \big \}$ indexed by more
general semigroups $\Phi \subseteq \RR^{\ge }$.  For
example, suppose that
$\Phi = \NN + \NN
\sqrt{2}$, and put 
\[ \fra_{j + k \sqrt{2}} \ = \ \frb^j \cdot \frc^k \]
for some fixed ideals $\frb, \frc \subseteq R$. Then
the volumes of the two $\NN$-graded subsystems $\{
\fra_{j} \}$ and
$\{ \fra_{k \sqrt{2}} \}$ will not in general coincide.
In order to avoid this sort of pathology, we will
henceforth adopt the following
\begin{convention}
For the remainder of this section, we deal with the
graded families arising from ideal filtrations; that is,  we
work with graded families $\big \{ \fra_{m} \big \}_{m
\in \Phi}$ which satisfy the additional condition
\begin{equation} \label{Filtration.Cond}
\fra_m \ \subseteq \ \fra_{m^\pr} \ \ \text{for any two
indices $m , m^\pr \in \Phi$ with }
\ \ m \ge  m^\pr. 
\end{equation}
\end{convention}
 \noi Of course (\ref{Filtration.Cond}) is automatic
for the graded families arising from valuations.  

\begin{remark}\label{finite} The volume of a graded
family of $\frmm$-primary ideals satisfying
(\ref{Filtration.Cond}) is  finite. Indeed, fix any $m
\in \Phi$. Since $
\fra_{m}^\ell \subset \fra_{m\ell} $ for every positive
integer $\ell$, we see that 
$$
\frac{\length(R/\fra_{m}^\ell)}{(m\ell)^n}\ \geq \
 \frac{\length(R/\fra_{m\ell})}{(m\ell)^n}
$$ for all $\ell \in \NN$. So taking the limit as $\ell$ 
gets large,
 Lemma \ref{vol.prop}  below  ensures that the volume of
$\{\fra_{\bullet}\}$  is bounded above by the rational
number $e(\fra_m)/m^n$, where $e(\fra_m)$ denotes the
multiplicity of the ideal $\fra_m$.  
\end{remark}

\begin{example}\label{vol.int}

\begin{itemize}
\item[(i).] The volume of the trivial graded family
$\{\fra^m\}$ of powers of a fixed ideal $\fra$ is equal
to the multiplicity of the ideal. Likewise, the volume
of $\{\overline{\fra^m}\}$ of integral closures of a
fixed power of an ideal $\fra$  is also the multiplicity
of $
\fra.
$
\item[(ii).]  Let $v$ be the valuation on $k(x, y)$ 
given by ''order of vanishing at the origin" 
 in Example \ref{val.ex}(i). As we have seen, the
valuation ideals of $\nu$  on $k[x, y]$ are given by
$\fra_m = (x, y)^m$.  In this case, the volume is the
multiplicity of the maximal ideal $(x, y)$  in $k[x, y]$,
 which is one. 
\item[(iii).] Let $v$ be the monomial valuation of
Example \ref{val.ex}(iii). Then each $\fra_m \subset
k[x, y]$ is  generated by  the monomials 
$x^a y^b,$
 where $a + \pi b \geq m$. So the length of $k[x,
y]/\fra_m$ is equal to the number of integer points in
the first quadrant of the Cartesian  plane  inside the
triangle bounded by the $a$-axis, the $b$-axis, and the
line $a + \pi b = m$. The area of this triangle is
roughly
$\frac{1}{2\pi} m^2$. Taking the limit,  the volume of
$\nu$ on $\spec k[x, y]$ is 
$\frac{1}{\pi}$. An evident modification of this shows
that any positive real number can occur as the volume
of a graded family of ideals in $k[x,y]$. 
\item[(iv).] Let $\nu$ be the arc valuation of Example
\ref{val.ex}(iv). Then the  quotients $R/q_m$  are
spanned by the residues of $1, x, x^2, \dots, x^{m-1}$,
so the length of  $R/q_m$ is $m$ and the volume of $\nu$
is zero. However, the $1$-volume is 1; Cf Remarks
\ref{p.vol} and \ref{p.vol.ex}.
\end{itemize}
\end{example}

\begin{remark} The volume of the family of ideals
associated to a divisorial valuation was considered by
Cutkosky and Srinivas in \cite{CS}. In the
two-dimensional case, they show that this invariant is
always a rational number: this essentially reflects the
existence of Zariski decompositions.
However Kuronya \cite{Kur} gives an example of a four
dimensional divisorial valuation with irrational
volume. His construction makes use of Cutkosky's curves
in $\PP^3$ having irrational Castelnuovo-Mumford
regularity \cite{C}. 
\end{remark}

\begin{remark}\label{p.vol} It is also possible to
define the $p$-volume of a graded system
$\fra_{\bullet}$ for any $p \leq n = \dim R$ as
$$
\limsup_{m\in \Phi} \,  
\frac{\length(R/\fra_m)}{m^p / p!}.
$$ In the current paper, we will not pursue this
further. See, however, Remark \ref{p.vol.ex}.
\end{remark}

\begin{lemma}\label{vol.prop} Let $\{\fra_t\}_{t\in
\Phi}$
 be a graded system of $\frmm$-primary ideals in a
Noetherian local ring $(R, \frmm)$ of dimension $n$
satisfying (\ref{Filtration.Cond}).  Then for any fixed
positive $m \in
\Phi$, 
$$
\limsup_{t \in \Phi} \frac{\length(R/\fra_t)}{t^n} \ = \
\limsup_{\ell \in \NN} \frac{\length(R/\fra_{m\ell})}
{(m\ell)^n}.
$$ In particular,  the volume of the $m$-th Veronese
graded subsystem
$\{\fra_{m\ell}\}_{\ell}$ is given by
$$ {\vol}\big(\{\fra_{m\ell}\}_{\ell \in \NN}\big )
\  =\  m^n \, {\vol  }\big(\{\fra_{t}\}_{t\in
\Phi}\big).
$$
\end{lemma}

\begin{remark} It follows from Lemma \ref{vol.prop} that
volume can be defined as 
$$
\limsup_{\ell \in \NN} \    
\frac{\length(R/\fra_{\ell m})}{(\ell m)^n / n!},
$$ where  $m$  is any fixed non-zero real number. In
particular, with the convention that $1 \in \Phi$,
 the volume is 
$$
\limsup_{m \in \NN}
\    \frac{\length(R/\fra_m)}{m^n / n!}.
$$
\end{remark}

\begin{proof}[Proof of Lemma \ref{vol.prop}]
For each 
$t \in \Phi$, we have
$$ m \cdot \left[ \frac{t}{m}\right ]
 \, \leq t \, < \, m \cdot \left[ \frac{t}{m}\right] + m,
$$ where $[\frac{t}{m}]$  denotes, as usual, the 
greatest integer less than or equal to the real number
$\frac{t}{m}$. Setting $\ell = [ \frac{t}{m}]$, we have
thanks to (\ref{Filtration.Cond}):
$$
\fra_{(\ell+1)m} \subset \fra_t \ \subseteq \ 
\fra_{\ell m},
$$ so that 
$$
\frac{\length(R/\fra_{(\ell + 1) m})}{t^n} \ \geq \ 
 \frac{\length(R/\fra_t)}{t^n} \ \geq \ 
\frac{\length(R/\fra_{\ell  m})}{t^n}. 
$$ Since $\lim_{t \ra \infty} \frac{t}{[\frac{t}{m}] m
}  = 1$, one has  
$$
\limsup_{\ell\rightarrow \infty}
 \frac{\length(R/\fra_{\ell  m})}{(\ell m)^n} \ = \ 
\limsup_{t \rightarrow \infty}
 \frac{\length(R/\fra_{\ell m})}{t^n},
$$ and likewise with $\ell$ replaced by $\ell + 1$. Thus 
$$
\limsup_{\ell\rightarrow \infty}
 \frac{\length(R/\fra_{(\ell + 1) m})} {((\ell + 1)
m)^n}\  \geq \ 
\limsup_{t\rightarrow \infty}
\frac{\length(R/\fra_t)}{t^n} \ \geq \ 
\limsup_{\ell\rightarrow \infty}
\frac{\length(R/\fra_{\ell  m})}{(\ell m) ^n}. 
$$ Since the expression on the left here is equal to the
expression on the right, the lemma is proved.
\end{proof}

The next proposition shows that from the point of view
of  multiplier ideals, graded families with zero volume
are  trivial.

\begin{proposition}\label{triv.mult} Let
$ \fra_{\bullet}  $ be a graded family
 of $\frmm$-primary ideals in a
local ring
$(R, \frmm)$ essentially of finite type over a field of
characteristic zero. Assume that $ \fra_{\bullet} $
 satisfies (\ref{Filtration.Cond}). If $\{
\fra_{\bullet} \}$ has volume zero, then each of its
asymptotic multiplier ideals
$\frj_m(\fra_{\bullet})$ is the unit ideal.  
\end{proposition}

\begin{proof} Fix any $m \in \Phi$. Then for all $\ell
\in \NN$, we have
$$
\fra_{m\ell} \ \subseteq \  \frj_m^\ell, 
$$ whence $e(\frj_m) \ge \vol( \fra_{\bullet}) $
thanks to the previous Lemma. The assertion follows. 
 \end{proof}

The following proposition, combined with the results of the
previous section, proves Corollary C from the
introduction. 

\begin{proposition} \label{key}
 Let $(R, \frmm)$ be a regular local ring of dimension
$n$, essentially of finite type over a field of
characteristic zero.  Let $\{\fra_m\}_{m \in \Phi}$  be
a  graded family of $\frmm$-primary ideals of
$R$ satisfying (\ref{Filtration.Cond}), 
 and let $\{\frj_m\}_m$ be the associated sequence of
asymptotic multiplier ideals. Assume that there is a
fixed  non-zero element $\delta \in R$ such that
\begin{equation} \label{key.lemma.eqn}
 \delta \cdot \frj_m \ \subseteq \ \fra_m \ \ \text{for
all \ } m \in \Phi. \end{equation} Then 
\[\vol ( \fra_{\bullet} ) \ = \ \limsup_{m \rightarrow
\infty}
\frac {e(\fra_m)}{m^n} \ = \ 
\limsup_{m \rightarrow \infty} \frac {e(\frj_m)}{m^n}, \]
  where  $e(\fra)$ denotes the multiplicity of the
ideal $\fra $ in the local ring $R$.
\end{proposition}
\noi  (As in Lemma
\ref{vol.prop}, the limits can be taken over all 
$m \in \Phi$ or just over all positive integer multiples
of a fixed element in $\Phi$. For graded systems indexed by the natural
numbers $\NN$, one does not need to assume the filtration condition
 (\ref{Filtration.Cond}).)

\begin{proof}  Given any index $m \in \Phi$, set   $\frd_m
= \big(  \fra_m : \frj_m  \big)$. This is an
$\frmm$-primary ideal, and since all the $\frd_m$
contain the fixed element $\delta$ one verifies that
\[
\limsup_{m \to \infty} \, \frac{e(\frd_m)}{m^n} \ = \
0. \]
Indeed, since $\fra_1^m \subset \fra_m \subset \frd_m$,
we see that $\big(\fra_1^m + \delta\big) \subset \frd_m$
for all $m$. Thus $e(\frd_m) \leq 
e(\overline{\fra_1}^m) = m^{n-1}e(\overline{\fra_1})$,
where $\overline{\fra_1}$ is the image of the ideal $\fra_1$ in the
 $(n-1)$-dimensional ring $R/(\delta)$.

Now fix a large index
$m \in
\Phi$. Then for all $\ell$  we have
$$ (\frd_m   \frj_m)^\ell \ \subseteq \  
\fra_m^\ell \ \subseteq \ \fra_{m\ell} \ \subseteq \
 \frj_{m\ell}
\ \subseteq \ \frj_{m}^\ell ,
$$ and so 
$$
\frac{\length(R/ \frj_m^\ell)}{\ell^n} 
\ \leq \ \frac{\length(R/\fra_{m\ell})}{\ell^n} \ \leq \
\frac{\length(R/ \fra_m^\ell)}{\ell^n}\  \leq\ 
\frac{\length(R/(\frd_m    \frj_m)^\ell)}{\ell^n}
$$ for all $\ell$. Taking the limit as $\ell $ goes to
infinity, we find that 
\begin{equation}\label{eq2} e( \frj_m)\  \leq \ {\text{
vol}} (\{\fra_{m\ell}\}_\ell)\  \leq \  e(\fra_m) \
\leq \  e(\frd_m  \frj_m),
\end{equation}
 where $\{\fra_{m\ell}\}_\ell$ is  the $m^{th}$ Veronese
subgraded sequence of
$\{\fra_{\bullet}\}$. 

But now note that
 dividing by $m^n$, the  expressions on the left and
right here  (namely $\frac{e( \frj_m)}{m^n} $ and
$\frac{e(\frd_m \frj_m)}{m^n}$)  have the same limit
superior as $m$  gets large.  In fact,   Teissier's
 Minkowski Inequality  \cite[p. 39]{Tes}  implies that
$$ e(\frd_m \frj_m)^{1/n} \ \le \ e(\frd_m)^{1/n} +
e(\frj_m)^{1/n}
$$ for all  $m$, whence
$$
\frac{e(\frd_m  \frj_m)^{1/n}}{m} \ \leq \
\frac{e(\frd_m)^{1/n}}{m}
 + \frac{e(\frj_m)^{1/n}}{m}
$$ for all $m$. But  
$\lim_{m \rightarrow \infty} \frac{e(\frd_m)}{m^{n}} =
0, $ and so 
 $$
\limsup_{m \rightarrow \infty}\frac{e(\frd_m
\frj_m)}{m^n} 
\ \leq \
 \limsup_{m \rightarrow \infty} \frac{e(\frj_m)}{m^n}.
$$ On the other hand,  since $e(\frj_m) \leq
e(\delta_m\frj_m)$, the reverse inequality always holds. 
 Finally, using  Lemma \ref{vol.prop},
 we  conclude from (\ref{eq2})  that 
$$
\limsup_{m\rightarrow \infty} \frac{e(\frj_m)}{m^n} \,\,
= \,\,   {\text{ vol}}(   \fra_\bullet   ) \  = \ 
 \limsup_{m\rightarrow \infty} \frac{e(\fra_m)}{m^n},
$$ as claimed.
\end{proof}

\smallskip

\begin{remark} The proof shows  also that the ``volume''
of $\{\frj_m\}$ is
 equal to the volume of the graded system $\{\fra_m\}$
 (even though $\{\frj_m\}$ itself is not a graded
system).
\end{remark}

\begin{remark} 
 Another example of a graded family satisfying (\ref{key.lemma.eqn})
  is given by the base
loci of the linear series of a big divisor. Specifically, fix a
big divisor  $L$ on a smooth projective  variety  $X$,
and let $\frb^\pr_m \subset \OO_X$ be the base ideal of
the linear system $|mL|$.  The components of the base
locus stabilize as $m \rightarrow \infty$, so choose
one component and localize along it to get a graded
family of ideals  $\frb_m$ in the local ring along the
generic point of the
 component. One can show that there exists $D$ such that
$\OO_X(-D)\cdot \frj_m \subset
\frb_m$ for all $m$ (see \cite{PAG}, Chapter 10).  So the
conclusion of  Proposition \ref{key} holds   for
$\frb_{\bullet}$.
\end{remark}

\begin{remark} \label{tightbound}
The conclusion of  Proposition \ref{key} holds under the weaker
assumption that there is a family of
non-zero $\frmm$-primary ideals $\frd_m$, with $\limsup
\frac{e(\frd_m)}{m^n} = 0$, such that $\frd_m \cdot
\frj_m \subseteq \fra_m$. As far as we know,
 it is possible that {\it every} graded system ${\fra_{m}}$
has this property, namely that the sequence $\fra_m$ is ``tightly bound''
to the sequence $\frak j_m$ in the sense that 
$\limsup \frac{e(a_m:j_m)}{m^n}$ tends to zero as $m$ goes to infinity. 
In particular, Proposition \ref{key}
may hold for a completely arbitrary graded system of $\frm$-primary ideals,
 so in particular, for an arbitrary valuation. 
Since we posed this question in an earlier version of this manuscript,
 Mircea Mustata has shed some light on the question of whether Proposition
\ref{key} holds more generally.
Specifically, he shows that the volume of an arbitrary graded system of
ideals is equal to the limit of the normalized multiplicities 
$\frac {e(\fra_m)}{m^n}$ in general---
that is, that  the first equality in the conclusion of Proposition \ref{key}
holds without the assumption that there exists a $\delta$ such that 
$\delta \frj_m \subset \fra_m$ for all $m$; see \cite{Must}. 
In particular, Mustata shows that
 Corollary C holds for any rank one valuation, 
Abhyankar or not. 
 However,  the relationship with the sequence
of asymptotic multiplier ideals (as well as the second equality in Proposition
\ref{key}) remains an  open question in the general case. 
\end{remark}

Theorem \ref{main} implies that the volume of an
Abhyankar valuation on any model is positive, and we saw
in Example \ref{vol.int} (iv) that the volume of a
non-Abhyankar valuation can be zero.  However, it is not
the case that a valuation has positive volume if and
only if it is Abhyankar, as the example below shows.

\begin{example}\label{vol0} In \cite[pp. 102--104]{ZS},
there is a  construction of  a  valuation on
$K = k(x, y)$ with value group an arbitrary additive
subgroup of the rational numbers; see also \cite[\S 10, Example 12]{Vaq}.
  Using this, we can construct a
non-Abhyankar  valuation of arbitrary volume (even
normalizing so that $\nu(\frmm) = 1$).

Let $\nu(y) = 1$ and set $\nu(x) = \beta_0 > 1$, some
rational  number. Let $c_0$ be the smallest positive
integer such that 
$c_0\beta_0 \in \ZZ$. As in \cite{ZS}, there exists a 
valuation so that the polynomial
$$ q_1 = x^{c_0} - y^{\beta_0c_0}
$$ has value $\beta_1$  equal to  any rational number
greater than (or equal to)  the 'expected value' of
$\beta_0c_0$. Let us choose this value so that $\beta_1
= \frac{d_1}{c_1} > \beta_0c_0$, where
$d_i$ and $c_i$ are relatively prime positive integers
with $c_1$ relatively prime to 
$c_0$.

This process can be repeated, so that we can construct a
valuation having the values on $x, y, $  and $ q_1$ as
already specified, and having arbitrary rational value
$\beta_2 \geq \beta_1c_1$ on the polynomial
$$ q_2 = q_1^{c_1} + y^{\beta_1c_1}.
$$ Again, we make this choice of $\beta_2$  so that the
smallest positive integer $c_2$ such that $c_2\beta_2
\in \ZZ$ is relatively prime each of the  preceding
$c_i$.

In this way, we inductively construct a sequence of
polynomials $q_i$,  rational numbers
 $\beta_i$, and positive integers $c_i$ with the
following properties:
$$ q_{i+1} = q_i^{c_i} + y^{\beta_ic_i},
$$ and 
$$
\beta_{i+1} > \beta_ic_i,
$$ where $c_i$
 is the smallest positive integer such that 
 $c_i\beta_{i} \in   \ZZ$, and $c_i$  is relatively
prime to  the product $c_0c_1 \dots c_{i-1}$.
 As shown (even more generally) in \cite{ZS}, this
uniquely defines a valuation  $\nu$ on 
$k(x, y)$, such that $\nu(q_i) =  \beta_i$. Indeed,
using the Euclidean algorithm,  and setting $q_{-1} = y
$ and 
$q_0 = x$, 
 every polynomial has a  unique expression  as a sum of
'monomials' in the $q_i$:
$$
 q_{-1}^{a_{-1}} q_0^{a_0} q_1^{a_1} \dots q_t^{a_t}
$$ where 
 $a_{-1}$ is arbitrary but the remaining exponents $a_i$
satisfy 
$a_i < c_i.$ Because each of these 'monomials' has a
distinct value,
 the valuation ideals on $k[x, y]$ have the form
$$
\fra_m = \left(\{q_{-1}^{a_{-1}} q_0^{a_0} q_1^{a_1}
\dots q_t^{a_t} \, | \,
\sum_{j=-1}^t \beta_i a_i \geq m; \,\,  a_j \leq c_j-1
{\text { for }}  j \geq 0
\}\right).
$$

In particular, the quotients $k[x,y]/\fra_m$  have
vector space basis consisting of 'monomials' 
$$ q_{-1}^{a_{-1}} q_0^{a_0} q_1^{a_1} \dots q_t^{a_t}
\,  \, {\text{  where }}
\, \, 
\sum_{j=-1}^t \beta_i a_i < m\,\, {\text{   and  }} \,
\,  
 a_j \leq c_j-1 {\text { for }}  j \geq 0.
$$ Although the number of products $t$  here can be
arbitrary,  note that for each fixed $m$, we only need
$t$ up to the greatest integer such that $\beta_{t} < 
m$.

 So computing the volume amounts to counting the number
of monomials in this basis. A computation shows that 
the volume of $\nu$  is  the limit of the following 
decreasing sequence of rational numbers
$$
\alpha_i = \frac{1}{\beta_0}
\left(\frac{c_0\beta_0}{\beta_1}\right)
\dots \left(\frac{c_i\beta_i}{\beta_{i+1}}\right).  
$$ (Calculation hint:  this is the limit of the 
  subsequence   $\{2! \frac{\lambda(k[x,
y]/\fra_m)}{m^2}\} $ indexed by  $m = c_t\beta_t$, which 
bounds the volume below. On the other hand, we can
approximate
$\nu$ by a sequence of valuations 
 $\nu_i$ which take the values $\beta_j$ on $q_j$
 for $j \leq i$, and the 'expected values' on the
remaining $q_i$  (which is to say that the corresponding
sequences of $\beta_j$ and $c_j$ satisfy  $\beta_j =
c_{j-1}\beta_{j-1}$ for $j > i$). Then just observe that
the volume of each $\nu_i$  bounds the volume of 
$\nu$ above and 
 compute that the volume of $\nu_i$ is $\alpha_i$.)

By choosing the values of $\beta_i$ and $c_i$
appropriately,  one can make this limit be  any
non-negative real number (note we've normalized so that
$\beta_{-1} = \nu(y) = 1$). For a completely explicit
example, take $c_i$  to be the standard
 enumeration of the prime numbers ($c_0 = 2, c_1 = 3,
\dots$), and set $\beta_{i+1} = c_i\beta_i +
\frac{1}{c_{i+1}}$. Then the volume turns out to be the
reciprocal of the  infinite sum
$$ 1 + \frac{1}{c_0} +  \frac{1}{c_0c_1} +
\frac{1}{c_0c_1c_3} + \dots,
$$ which is approaches a real number between
$\frac{1}{2}$ and 1.
\end{example}

\begin{discussion}[The Associated Graded Algebra of
a Valuation]
 Fix a rank one valuation $\nu$  on a function field
$K/k$
 centered a local  domain $(R, \frmm)$ and let $\Phi$ be
the corresponding value semigroup $\nu(R)$. The {\it
associated graded algebra\/} of the valuation on $R$  is 
$$ {\text{gr}}_{\nu} R = \bigoplus_{m\in \Phi}
\fra_{m}/\fra_{>m},
$$ where $\fra_{>m}$ denotes the valuation ideal 
$\{f \, |  \, \nu(f) > m\}.$ It is easy to check that
${\text{gr}}_{\nu} R$ is a domain, but it is not
finitely generated over $R/\frmm$ (its  degree zero
piece)  in general. The transcendence degree of 
${\text{gr}}_{\nu} R$ over 
$R/\frmm$  is equal to the rational rank of $\nu$ plus
the  transcendence degree of $R_{\nu}/\frmm_{\nu}$ over
$R/\frmm$, where $\frmm_{\nu}$ denotes the maximal ideal
of the valuation ring $R_{\nu}$. If ${\text{gr}}_{\nu} R
$ is finitely generated, therefore, its Krull dimension
is equal to this sum. In this case, the associated
graded ring has dimension 
 equal to the dimension of $R$ if and only if the
valuation is Abhyankar.{\footnote{Here we are using a
slightly different, but equivalent (for function fields),
 form of Abhyankar's inequality (\ref{ab})  which says
that  $rat.rank \nu + trans.deg_{R/\frmm}
(R_{\nu}/\frmm_{\nu})
 \leq R$; this inequality differs from (\ref{ab})  by
addition of the same number, namely the transcendence
degree of  
 $R/\frmm$
 over $k$,
 to both sides.}}

When the associated graded ring  of $\nu$  is finitely
generated and 
$\NN$-graded, the volume of $\nu$ has a simple interpretation in
terms of ${\text{gr}}_{\nu} R$, namely it is  
 equal to the Hilbert
multiplicity of $\nu$. Recall that if $A = \oplus A_m$ is a finitely
generated $\NN$-graded domain  over a field $A_0$ containing
non-zero elements of every degree, then
there exists a positive rational number $e$ such that
\begin{equation}
\dim A_m = e \frac{m^{n-1}}{(n-1)!} + O(m^{n-2})
\end{equation} where $n$ is the Krull dimension of $A$
and $\dim$ denotes the dimension over the field $A_0$.
This number $e$  is called the {\it Hilbert
multiplicity} of 
$A$.  [To see this, note that for some $r$, the Veronese subalgebra
$A^{(r)}$ is generated in degree one (or $r$), so $A$ decomposes
as a direct sum of $A^{(r)}$-modules 
$A_{(0)} \oplus A_{(1)} \oplus \dots \oplus A_{(r-1)},$ 
where $A_{(i)} = \oplus_{j\in \NN} A_{jr +  i}$.
Thus each of the modules $A_{(i)}$ has some multiplicity $e_i$ over the ring 
$A^{(r)}$. In general, these multiplicities can be different,
but if $A$ has elements of every degree, one shows that the $e_i$ are all
 equal. Indeed, to see that $e_i = e_{i'}$, just take a non-zero 
 element $x \in A_{i-i'}$ and note that the  cokernel of the injective
map $A_{(i')} \overset{x}{\rightarrow} A_{(i)}$ has dimension strictly less
than the common dimension of  $A_{(i)}$ and $A_{(i')}$. Thus 
 $A_{(i)}$ and $A_{(i')}$ necessarily have the same multiplicity, $e$.
Thus the Hilbert polynomials of all the $A_{(i)}$ have the same leading 
terms, leading to the formula for the dimension above (with $e$ suitably
normalized).]

 Thus
$$
\length (R/\fra_m) \, = \, \sum_{i=0}^{m-1}
\length(\fra_m/\fra_{m+1}) \,  = \, \frac{e}{(n-1)!}
\left[\sum_{i=0}^{m-1} i^{(n-1)}  + O(m^{n-2})\right],
$$
 and since $\sum_{i=0}^{m-1}[ i^{(n-1)} + O(m^{n-2})] =
\frac{m^n}{n} + O(m^{n-1})$, we see that 
$$
\lim_{m\ra\infty}\,   \frac{\length(R/\fra_m)}{m^n / n!}
= e.
$$
\begin{remark}\label{p.vol.ex} When ${\text{gr}}_{\nu} R
$ is finitely generated, the preceding discussion
indicates that it is natural to consider the $p$-volume
of $\nu$ on $R$, where $p = rat.rk. \nu +
trans.deg._{R/\frmm} (R_{\nu}/\frmm_{\nu})$; see Example
\ref{vol.int} (iv). However, Example \ref{vol0} indicates
that the $p$-volume need not be finite in general, even
when the graded ring has transcendence degree $p$. 
\end{remark}
\end{discussion}

\section{Generalizations and Further Applications}

Let $D$ be an arbitrary effective divisor on a smooth
variety $X$ and let $\pi: X \lra Y$ be a proper
birational map to a smooth variety $Y$ that  collapses
$D$ to a point. The collection 
\begin{equation}
\label{div}
\{\pi_*\OO_X(-mD) \}_{m\in \NN}
\end{equation} forms a graded family of ideals in
$\OO_Y$. Although this is not  the graded family of a 
valuation, it can  be handled by the  methods  developed
here for families of valuation ideals because it is an
intersection of graded families of valuation ideals.

\begin{definition} Let 
 $\left\{\{ a^{\lambda}_{m} \}_{m\in
\Phi}\right\}_{\lambda \in \Lambda}$
 be an arbitrary collection of graded families, all
indexed by 
 the same semigroup $\Phi \subset \RR$. The intersection
graded family is defined by 
$$
\bigcap_{\lambda \in \Lambda}  \fra_{\bullet}^{\lambda}
\,\, := \,\, 
\left\{\bigcap_{\lambda\in \Lambda}\fra_m^{\lambda}
\right\}_{m \in \Phi}.
$$
\end{definition}
\noi Note that if each $\fra^\lambda_{\bullet}$
satisfies (\ref{Filtration.Cond}), then so too does their
intersection.

 The asymptotic multiplier ideals of an intersection
family satisfy
\begin{equation}
\label{intersect}
\frj_m\left(\bigcap_{\lambda\in
\Lambda}\fra^{\lambda}_{\bullet}\right)
\subset 
\bigcap_{\lambda\in
\Lambda}\frj_m\big(\{\fra^{\lambda}_\bullet\}\big) 
\end{equation} Indeed, for each $\lambda$, we have 
$ \{\bigcap_{\lambda\in \Lambda}\fra_m^{\lambda} \}
\subset \fra_m^\lambda
$ for all $m$, so 
 the corresponding inclusion holds for the asymptotic
multiplier ideals.

 Let $\mathcal  {S}$  be any collection of rank one
valuations centered on a domain $R$. For each $m \in
\RR_{\geq 0}$, set  
\begin{equation}
\label{gf.coll}
\fra_m = \{f \in R \, | \, \nu(f) \geq m
\, \, {\text { for all }} \nu \in \mathcal  {S}\}.
\end{equation} The collection 
$\{\fra_m\}_{m\in \RR_{\geq 0}}$ forms a graded family
of ideals in $R$ indexed by the non-negative real
numbers  (or by considering only distinct such ideals,
some subsemigroup of $\RR$.)  This graded family  is the
intersection, over all $\nu \in \mathcal S$,  of the
graded families $\{q_m(\nu)\}_{m\in \RR_{\geq 0}}$ of
$\nu$-valuation
 ideals.

As a variant, one can also assign multiplicities to the
valuations in
$\mathcal  {S}$. Say for each $\nu$  in $\mathcal  {S}$,
we assign some positive real number  $e_{\nu}$. Then for
each $m \in \RR$, set
\begin{equation}
\label{gf.coll.mult}
\fra_m = \{f \in R \, |  \,\,\, \frac{\nu(f)}{e_{\nu}}
 \geq  m \, \,\,
 {\text { for all }} \nu \in \mathcal  {S}\}.
\end{equation} The collection 
$\{\fra_m\}_{m\in \RR}$ forms a graded family of ideals
in $R$, indexed by (some subsemigroup of) the real
numbers. For example,  the graded family (\ref{div})
above is of this form: If 
 $D = \sum e_i D_i$ where the $D_i$ are prime divisors,
then the set $\mathcal S$  is the set of valuations
$\nu_i$  given by order of vanishing along $D_i$ and the
multiplicities $e_{\nu_i} $ are the coefficients $e_i$.
These graded families are also intersections:  the
graded family (\ref{gf.coll.mult})  is the 
intersection, over all $\nu$ in $\mathcal S$,  of the
graded families $\{q_{e_{\nu}m}(\nu)\}_{m}$ of the
$e_{\nu}$-th Veronese subfamily of the graded family of
valuation ideals of $\nu$.  Alternatively, it can be
interpreted as an  intersection of graded families of 
valuation ideals:  it is the intersection of the graded
families of 
 valuation ideals in the set $\mathcal S'$, where the
set $\mathcal S'$ is obtained from $\mathcal S$ by
replacing each valuation $\nu$  in $\mathcal S$
 by the valuation
 $\frac{1}{e_{\nu}} \nu$. So this is really no more
general than the intersection (\ref{gf.coll})  discussed
in the previous paragraph. Once again, condition 
(\ref{Filtration.Cond}) is satisfied by these families. 

\medskip One can also define the product of two 
 graded families $\{\fra_{\bullet}\}_{m \in \Phi}$
 and  $\{\frb_{\bullet}\}_{m\in \Phi}$
  by  
$$
\{\fra_m\frb_m\}_{m\in \Phi}.
$$ Note that the asymptotic multiplier ideals of the
product graded family satisfy
\begin{equation}
\label{product}
\frj_m(\fra_{\bullet}  \frb_{\bullet}) \ \subseteq  \
\frj_m{(\fra_{\bullet})} \frj_m{( \frb_{\bullet})},
\end{equation}
 since,  for all large $p$,  we have the inclusion of
'usual' multiplier  ideals
$\MI{\frac{1}{p}\cdot (\fra_{mp}\frb_{mp}) }\subset 
\MI{\frac{1}{p}\cdot \fra_{mp} }\MI{\frac{1}{p}\cdot
\frb_{mp}} $ by subadditivity \cite{DEL}. By induction,
the product of any finite number of graded families
indexed by the same semi-group is defined, and the same
multiplicativity of the asymptotic multiplier ideals
holds.

Our Main Theorem \ref{main} and its ensuing corollaries
 can be extended to graded families arising as  finite
intersections or products
 of Abhyankar valuations.  Explicitly:

\begin{corollary}\label{gen} Let $\mathcal{S}$ be a
finite  collection   of rank one Abhyankar  valuations 
of a function field $K/k$ of characteristic zero, all
centered on some local ring $R$ of a model of $K/k$.
 Let $\{\fra_{\bullet}\}$ denote either the 
corresponding  intersection or the corresponding  
product graded  family of ideals.  Then  there exists a
non-zero element $\delta \in R$ such that 
$$
\delta \, \frj_m \  \subset  \ \fra_m
$$ for all $m \in \RR$, where $\{\frj_m\}$ is the
associated sequence of asymptotic multiplier ideals.
Furthermore, the conclusions of Theorem A  and Corollary C 
  hold for this graded
family of ideals.
\end{corollary}

\begin{proof} This is immediate from Theorem \ref{main}.
For each $\nu \in \mathcal{S}$, let 
$\delta_{\nu}$  be the non-zero element guaranteed by
Theorem \ref{main} and let $\delta$ be their product.
 The desired inclusion follows from
(\ref{intersect}) or from (\ref{product}), respectively,
for the intersection and the product case.
\end{proof}

\begin{example} For the graded system (\ref{div}) of a
divisor $D$, Corollary \ref{gen} immediately implies Corollary B from the
Introduction. It also guarantees the
existence of a non-zero $\delta$ such that 
$$
\delta^{\ell}  (\pi_*\OO_X(-m\ell D)) \ \subseteq \
(\pi_*\OO_X(-mD))^{\ell} 
$$ $\ell, m \in \NN$ with $m\gg 0 $.
\end{example}

Finally,  we point out the following Minkowski
Inequality for intersections  and products of graded
families.

\begin{corollary}[Minkowski Inequality]
\label{mink} Let $\fra_{\bullet}$ and
$\frb_{\bullet}$ be two graded families of
$\frmm$-primary ideals in a regular local  ring
essentially of finite type over a field of
characteristic zero. Assume that the hypothesis of Proposition
\ref{key} holds for each of these systems. Then
$$ {\textnormal{vol}} (\fra_{\bullet} \bigcap
\frb_{\bullet})^{1/n} \ \leq\ 
 {\textnormal{vol}} (\fra_{\bullet} 
\frb_{\bullet})^{1/n}
\ \leq \  {\textnormal{vol}} (\fra_{\bullet})^{1/n} +
{\textnormal{vol}} (\frb_{\bullet})^{1/n}
$$
\end{corollary}

\begin{proof} Since $\fra_m \frb_m \subseteq \fra_m
\bigcap \frb_m$ for all $m$, it is evident that 
$ {\text{vol}} (\fra_{\bullet} \bigcap \frb_{\bullet})
\ \leq \ 
 {\text{vol}} (\fra_{\bullet}  \frb_{\bullet}).
$ So it is enough to prove the inequality for product
families.

It follows from the subadditivity relation
(\ref{product})
 that if the hypothesis of Proposition
\ref{key} holds for each of $\fra_{\bullet}$ and
$\frb_{\bullet}$ then it holds also for the product
system. This being said, the result follows from 
Proposition \ref{key} and Teissier's Minkowski inequality for
multiplicities. Indeed,
\begin{align*} {\text{vol}}
(\fra_{\bullet}\frb_{\bullet})^{1/n} \ &= \ 
\left[\limsup_m \frac{e(\fra_m\frb_m)}{m^n}
\right]^{1/n}  \\   &\leq \
\limsup_m \left\{ \left[
\frac{e(\fra_m)}{m^n}\right]^{1/n}
 + \left[\frac{e(\frb_m)}{m^n}\right]^{1/n} \right\} \\
&\leq \
\limsup_m \left\{ \left[
\frac{e(\fra_m)}{m^n}\right]^{1/n}\right\} \, 
 + \, \limsup_m
\left\{ \left[\frac{e(\frb_m)}{m^n}\right]^{1/n}
\right\}
\\ &=  \  {\text{vol}}(\fra_{\bullet})^{1/n } +  {\text{
vol}}(\frb_{\bullet})^{1/n},  
\end{align*} as required.
\end{proof}

 \begin{remark} Note that the hypothesis of Proposition
\ref{key} is satisfied in particular for   graded
families that are finite
 products or intersections of graded families of
Abhyankar valuations. This is the content of Corollary
\ref{gen}.  More generally, the product or intersection
graded family of  any finite collection of graded
families that satisfy the hypothesis also
satisfy the hypothesis.
 Indeed, we verified this 
in the proof of Corollary \ref{mink}  for products using
(\ref{product}) and the same argument, using
(\ref{intersect}) instead of
 (\ref{product}), works for intersections.
\end{remark}

\begin{example} As a special case of the Minkowski
inequality,
 fix an  $\frmm$-primary ideal  $\fra$  in a  local ring
$R$ of a smooth complex variety  and  let $\nu_1, \dots,
\nu_r$ be the associated Rees valuations of $\fra$. (The
definition of Rees valuations is recalled in the
Introduction.)  Set  $e_i = \min_{f \in I} \nu_i(f)$.
Then the associated graded family
$$
\fra_m =
\big \{f \in R \, \mid \, \frac{\nu_i(f)}{e_i}
 \geq m
\,  \, {\text { for }}  i = 1, 2, \dots, r \big \}
$$ is nothing more than the  graded family of integral
closures of powers of
$\fra$ of Example \ref{first.examples} (ii).  In
particular, its volume is the multiplicity of $\fra$.
This family is a finite intersection of graded families
of Abhyankar valuations, so Corollary \ref{mink}  can be
applied. This produces the bound
$$ e(\fra)^{1/n} \ \leq \ e_1 {\text{ vol}}(\nu_1)^{1/n}
+
\dots +   e_r {\text{  vol}}(\nu_r)^{1/n}.
$$
\end{example}

\end{document}